\documentclass{amsart}
\usepackage{color}

\usepackage{amsmath} 
\usepackage{amssymb}
\usepackage{mathrsfs}

\newtheorem{theorem}{Theorem}[section] 
\newtheorem{claim}[theorem]{Claim}
\newtheorem{mc}[theorem]{Main Claim}

\theoremstyle{definition}
\newtheorem{definition}[theorem]{Definition}
\newtheorem{example}[theorem]{Example}
\newtheorem{observation}[theorem]{Observation} 

\newtheorem{discussion}[theorem]{Discussion}

\newtheorem{question}[theorem]{Question}

\newtheorem{hypothesis}[theorem]{Hypothesis}

\theoremstyle{remark}
\newtheorem{remark}[theorem]{Remark}
\newtheorem{notation}[theorem]{Notation}

\newtheorem{context}[theorem]{Context}

\newcommand{\rest}{{\restriction}}
 
\newcommand{\inc}{{\rm inc}} 
\newcommand{\tp}{{\rm tp}} 

\newcommand{\Rang}{{\rm Rang}}
\newcommand{\set}{{\rm set}}
\newcommand{\wilog}{{\rm without loss of generality}}
\newcommand{\Wilog}{{\rm Without loss of generality}}

\newcommand{\then}{{\underline{then}}}
\newcommand{\when}{{\underline{when}}}

\newcommand{\mn}{{\medskip\noindent}}
\newcommand{\sn}{{\smallskip\noindent}}

\newcommand{\bbB}{{\mathbb B}}

\newcommand{\gC}{{\mathfrak C}}

\newcommand{\cH}{{\mathcal H}}

\newcommand{\bbL}{{\mathbb L}}

\newcommand{\cP}{{\mathcal P}}

\newcommand{\bbQ}{{\mathbb Q}}

\newcommand{\varp}{{\varepsilon}}
 
\newcommand{\cU}{{\mathcal U}}

\newcommand{\cf}{{\rm cf}}

\newcount\skewfactor
\def\mathunderaccent#1#2 {\let\theaccent#1\skewfactor#2
\mathpalette\putaccentunder}
\def\putaccentunder#1#2{\oalign{$#1#2$\crcr\hidewidth
\vbox to.2ex{\hbox{$#1\skew\skewfactor\theaccent{}$}\vss}\hidewidth}}

\newenvironment{PROOF}[2][\proofname.]
   {\begin{proof}[#1]\renewcommand{\qedsymbol}{\textsquare$_{\rm #2}$}}
   {\end{proof}}

\usepackage[hidelinks]{hyperref}

\begin{document}
\makeatletter\def\shfiuwefootnote{\gdef\@thefnmark{}\@footnotetext}\makeatother\shfiuwefootnote{Version 2017-06-09\_12. See \url{https://shelah.logic.at/papers/886/} for possible updates.}

\title {Definable groups for Dependent and 2-dependent theories \\
Sh886}
\author {Saharon Shelah}
\address{Einstein Institute of Mathematics\\
Edmond J. Safra Campus, Givat Ram\\
The Hebrew University of Jerusalem\\
Jerusalem, 91904, Israel\\
 and \\
 Department of Mathematics\\
 Hill Center - Busch Campus \\ 
 Rutgers, The State University of New Jersey \\
 110 Frelinghuysen Road \\
 Piscataway, NJ 08854-8019 USA}
\email{shelah@math.huji.ac.il}
\urladdr{http://shelah.logic.at}
\thanks{The author would like to thank the Israel Science Foundation for
partial support of this research (Grant No. 242/03).
I would like to thank Alice Leonhardt for the beautiful typing. 
 First Typed - 05/Dec/15}

  % Previous version - Sept.20, 2016

  % Formerly F753 - changed July/06

\subjclass{Primary 03E17; Secondary: 03E05, 03E50}

\keywords {Model theory, groups, classification theory, dependent theories}

\date{June 9, 2017}

\begin{abstract}
Let $T$ be a (first order complete) dependent theory, ${\gC}$ a
 $\bar\kappa$-saturated model of $T$ and $G$ a definable subgroup 
which is abelian.  Among subgroups of bounded index which are the
 union of $< \bar\kappa$ type-definable subsets
 there is a minimal one, i.e. their intersection has bounded
index.  See history in \cite{Sh:876}.  
We then deal with definable groups for 2-dependent theories, a wider class
 of first order theories proving that for many pairs $(M,N)$ of models,
the minimal bounded subgroup definable over $M \cup N$ is the
intersection of the minimal ones for $M$ and for $N$.
\end{abstract}

\maketitle
\numberwithin{equation}{section}
\setcounter{section}{-1}
\newpage

\section {Introduction}

Assume that $T$ is a dependent (complete first order) theory, 
${\gC}$ is a $\bar \kappa$-saturated model of $T$ (a monster), $G$ is a
type-definable (in ${\gC}$) group in ${\gC}$ (of course we
consider only types of cardinality $< \bar \kappa$).

A subgroup $H$ of $G$ is called bounded if the index
$(G:H)$ is $< \bar \kappa$.
By \cite{Sh:876},
we know that among $\mathbf H_{\text{tb}} = \{H:H$ is a
type-definable subgroup of $G$ of bounded, i.e. with index $< \bar
\kappa\}$  there is a minimal one.

But what occurs for $\mathbf H_{\text{stb}} = \{H:H$ is a 
union-type-definable (see below) subgroup of $G$ of bounded, i.e. $< \bar
\kappa$ index$\}$? is there a minimal one?  Our main result is a
partial positive answer: if $G$ is an abelian group then in $\mathbf
H_{\text{stb}}$ there is a minimal one.

We call $C \subseteq {\gC}$ union-type-definable (over $A$) when
for some sequence $\langle p_i(x):i < \alpha\rangle$ we have $\alpha <
\bar\kappa$, each $p_i(x)$ is a type (over $A$; of course $A$ is of
cardinality $< \bar\kappa$) and $C = \cup\{p_i({\gC}):i < \alpha\}$
where $p_i({\gC}) = \{b \in {\gC}:b$ realizes $p_i({\gC})\}$;
this is equivalent to being $\bbL_{\infty,\kappa}$-definable in 
${\gC}$ for some $\kappa \le \bar\kappa$. 

In Definition \ref{dt.22} we recall the definition of 2-dependent
$T$ (where $T$ is 2-independent when some $\langle \varphi(\bar x,\bar
b_m,\bar c_n):m,n < \omega\rangle$ is an independent sequence of formulas); see
\cite[\S5 (H)]{Sh:863}.  Though a reasonable definition, can we say
anything interesting on it?  Well, we prove the following result related to
\cite{Sh:876}.

Let $G_A$ be the minimal type-definable over $A$ subgroup of $G$,
for suitable $\kappa$; as we fix $A$ it always exists.  Theorem
\ref{dt.35} says that if $M$ is
$\kappa$-saturated and $|B| < \kappa$ then $G_{M \cup B}$ can be
represented as $G_M \cap G_{A \cup B}$ for some $A \subseteq M$ of
cardinality $< \kappa$.  So though this does not prove ``2-dependent
is a dividing line", it seems enough for showing it is an interesting property.

The first theorem on this line for $T$ stable is
of Baldwin-Saxl \cite{BaSx76}.
Recently Hrushovski, Peterzil and Pillay \cite{HPP0x} investigated
definable groups, o-minimality and measures; an earlier work on
dependent theories is \cite{Sh:715}; and
on definable subgroups in o-minimal $T$ is Berarducci, Otero, Peterzil
and Pillay \cite{BOPP05} where the existence of
the minimal type-definable bounded
index theorem and more results are proved for o-minimal theories.

A natural question was whether there is a minimal type-definable 
bounded subgroup,
when $T$ is dependent. Assuming more (existence of measure) on the
group, related suitably to the family of definable sets, this was proved
in the original version of Hrushovski-Peterzil-Pillay \cite{HPP0x}.  
Then \cite{Sh:876} proves this for every dependent theory and definable
group.
(The final version of their paper \cite{HPP0x} includes 
an exposition of the proof of \cite{Sh:876}.)

Recent works of the author on dependent theories are \cite{Sh:783}
(see \S3,\S4 on groups) \cite{Sh:863} (e.g. the first order theory of
the $p$-adics is strongly$^1$ dependent but not strongly$^2$ dependent,
see end of \S1; on strongly$^2$ dependent fields see \S5) and
\cite{Sh:900} and later \cite{Sh:950}.  

Hrushovski has pointed out the following application of \S2
to Cherlin-Hrushovski [ChHr03].
They deal with complete first order theories $T$ with few
$\lambda$-complete 4-types such that every finite subset of $T$ has a
finite model.  Now from the classification of finite simple  groups 
several (first order) properties of such theories were deduced.
\relax From them they get back information on automorphism groups of finite
structures.
An interesting gain of this investigation is that if looked at this 
as a round trip
from classification of finite simple groups to such $T$'s, we get
uniform bounds for some of the existence results in the classification
of finite simple groups.
Note that they deal with the case the finite field is fixed, while the
vector space over it varies.
This is related to  the properties being 
preserved by reducts (and interpretations) of
first order theories hence we get the uniform bounds.
Now by the present work some of those first 
order properties of such $T$'s which they proved using 
finite simple groups, 
are redundant in the sense that they follows from the others
though originally it was non-trivial to prove them.  Why?  Because (the
relevant first order theories are 2-dependent by Observation
\ref{nd.17} and so) Theorem \ref{dt.35} can be applied.  More
specifically, the main point is a version of modularity that looking at
an expansion of $M = (V,B)$, where $V$ is a vector space over a fixed
finite field $F,B$ a bilinear map from $V$ to $F,M$ an ultraproduct of
finite structures, for suitable $M_0 \prec M_\ell \prec M$ for
$\ell=1,2,acl(M_1 \cup M_2) = \{a_1 + a_2:a_1 \in M_1,a_2 \in M_2\}$.  

This raises further questions.

We thank Itay Kaplan for pointing out many many things which should be
corrected and the referee for suggestions to make the paper shorter and more
user-friendly.

The reader may ask

\begin{question}
Is the notion of 2-independent and more generally
$n$-independent first order theory interesting?  

In my view they certainly are.  Why?
\medskip

\noindent
(A) \quad Their definition is simple and natural.  

True this is a statement easy to argue about, hard to convince about, 
you may say it is a matter of personal taste. But 
we can look at history: earlier cases in which the
author suggest such notions proving something on them.  True, the
interest in the (failure of the) strict order property (see \cite{Sh:12};
\cite[\S5]{Sh:702}, Kikyo-Shelah \cite{Sh:748}, Shelah-Usvyatsov, 
\cite{Sh:789}) and $n$-NSOP relatives (see \cite[\S2]{Sh:500}, 
Dzamonja-Shelah \cite{Sh:692}, Shelah-Usvyatsov \cite{Sh:844})
have not so far been proved to skeptics.  But other cases looked 
for a non-trivial
number of years to be in a similar situation but eventually become quite
popular: stable and superstable (\cite{Sh:2}), the independence
property itself (\cite{Sh:12}), and a more strict version of the
(failure of the) strict order property = simple theory (= failure of
the tree property).
\medskip

\noindent

(B) \quad The theory of such classes is not empty; we prove in \S2
that ``the minimal bounded subgroup of $G$ over $A$" behaves nicely.
A skeptical reader suggests that to justify our interest we have to
point out ``natural examples", i.e. from other branches of mathematics.  We
object to this, just as we object to a parallel view on mathematics as a whole.
In both cases applications are a strong argument for and desirable,
\underline{but} not a necessary condition.  In fact such cases were not
present in the author first works on the notions mentioned above.

Moreover the connection to Chelin-Hrushovski \cite{ChHr03} already
showed such connection.

In any case there are such examples, 
generally if in ${\gC} = {\gC}_T$ we can define a
field $K$ and vector spaces $V_1,V_2,V_3$ over it (or just abelian
groups) and a bi-linear mapping $F:V_1 \times V_2 \rightarrow V_3$,
then in the non-degenerated case the formulas $\{F(x,b) = 0_{V_3}:b
\in V_2\}$ form an infinite family of independent formulas.  As
bi-linear maps appear in many examples, it makes sense even to a
reasonable skeptic to look for a (real) parallel to the class of
stable theories including them and 2-dependent is the most natural one.
\end{question}

\centerline {$ * \qquad * \qquad *$}

\noindent
Concerning \cite{Sh:876}, note that
\begin{claim}
\label{2n.7}  
[$T$ dependent]  

Assume
\mn
\begin{enumerate}
\item[$(a)$]   $G$ is a $A^*$-definable semi-group with cancellation
\sn
\item[$(b)$]    $q(x,\bar a)$ is a type, $q({\gC},\bar a)$ a
sub-semi-group of $G$
\sn
\item[$(c)$]    for $\bar a'$ realizing $\tp(\bar a,A^*)$ let
  $\set(\bar a') = \{b$: for some $a,c \in q(\gC,\bar a')$ we have $a
  b c \in q(\gC,\bar a')\}$; (note that if $G$ is a group, $q(\gC,\bar a)$ a
  sub-group then set$(\bar a') = q(\gC,\bar a))$
\end{enumerate}
\mn
\then \, we can find $q^*$ and $\langle \bar a_i:i < \alpha \rangle$
and $B$ such that:
\mn
\begin{enumerate}
\item[$(\alpha)$]   $\alpha < \lambda^+$ where $\lambda = |T| + |\ell
  g(\bar a)|$
\sn
\item[$(\beta)$]   {\rm tp}$(\bar a_i,A^*) = \text{\rm tp}(\bar
a,A^*)$ hence $q({\gC},\bar a_i)$ a sub-semi-group of $G$
\sn
\item[$(\gamma)$]   $q^*(x) = \cup\{q(x,\bar a_i):i < \alpha\}$
\sn
\item[$(\delta)$]    $B \subseteq \bigcup\limits_{i < \alpha} q(G,\bar a_i)
 \subseteq G$ and $|B| \le |\alpha|$; in fact, $B \subseteq
 \bigcup\limits_{i} \, \bigcap\limits_{j \in \alpha \backslash \{i\}}
 q(G,\bar a_i)$ 
\sn
\item[$(\varepsilon)$]   if $\bar a'$ realizes 
{\rm tp}$(\bar a,A^*)$ and $B \subseteq q({\gC},\bar a')$ then
$q^*({\gC}) \subseteq \set(\bar a')$.
\end{enumerate}
\end{claim}

\begin{remark}
This was part of \cite{Sh:876}, but it was claimed
that it seems to be wrong: 
Note that we do not require ``$B \subseteq q^*({\gC})$".  The
following example may clarify the claim.

Let $T$ be the theory of ($\bbQ,+,<,0,1)$ say,
and $G$ the monster model.  Let $a$ be infinitesimal, namely it realizes
$p(x) := \{0 < x <1/n:n=1,2,3,\ldots\}$.  Let $q(x,a)$ say  that $|x| <
a/n$ for $n=1,2,\ldots$.  Then for any set $A \subseteq p({\gC})$
(of cardinality $< \bar\kappa$) and set $B$ (of cardinality $< \bar\kappa$)
 contained in $q({\gC},a')$ for all $a' \in A$, 
there is by compactness some $d$ realizing $p(x)$ such that 
$d/n >b$ for all $b \in B$ and
$n=1,2,\ldots$ and $d < a'/n$ for all $a' \in A$ and $n=1,2,\ldots$.
However, what $(\varepsilon)$ of \ref{2n.7} says is: if $a'$
realizes $p(x)$ (i.e. $a'$ infinitesimal) and $B \subseteq
q({\gC},a')$, i.e. $b \in B \wedge n < \omega \Rightarrow |b| < a'/n$ then
$q^*({\gC}) \subseteq q({\gC},a')$ so if $a \in B$ this holds.
\end{remark}

\begin{PROOF}{\ref{2n.7}} 
We try to choose $\bar a_\alpha,b_\alpha$ by induction
on $\alpha < (|T|^{\aleph_0})^+$ such that:
\mn
\begin{enumerate}
\item[$\circledast$]   $(a) \quad \bar a_\alpha$ realizes tp$(\bar a,A^*)$
\sn  
\item[${{}}$]     $(b) \quad b_\alpha \notin q({\gC},\bar a_\alpha)$ 
and $b_\alpha \in G$ and moreover $a,c \in q(\gC,a_\alpha)
\Rightarrow$

\hskip25pt $a b_\alpha c \notin q(\gC,a_\alpha)$, i.e. $b_\alpha
\notin \set(a_\alpha)$
\sn
\item[${{}}$]    $(c) \quad b_\alpha$ realizes $q(x,\bar a_\beta)$
for $\beta < \alpha$
\sn
\item[${{}}$]    $(d) \quad b_\beta$ realizes $q(x,\bar a_\alpha)$ 
for $\beta < \alpha$.
\end{enumerate}
\mn
First, if we are stuck at some $\alpha < (|T|^{\aleph_0})^+$ 
then the desired result is
exemplified by $\langle \bar a_i:i < \alpha\rangle$ and 
$B := \{b_i:i < \alpha\}$.  Note that if $\alpha \ge 1$ then $B
\subseteq \cup\{\cap\{q({\gC},\bar a_i):i \in \alpha \backslash
\{j\}:j < \alpha\}$ and if $\alpha=0$ then $B = \emptyset$.

Second, if we succeed we get contradiction similarly to the proof in
\cite[\S1]{Sh:876} but we elaborate; let $\bar y = \langle y_i:i <
\ell g(\bar a)\rangle$.  Let $\varphi_\alpha(x,\bar y) \in q(x,\bar
a_\alpha)$ be such that $\gC \models \neg \varphi_\alpha(b_\alpha,\bar
a_\alpha)$, exists by clause (b) of $\circledast$.  The number of
possible $\varphi_\alpha(x,\bar a_\alpha)$ is $\le \lambda$, hence for
some $\varphi_*(x,\bar y) \in \bbL(\tau_T)$ the set $W := \{\alpha <
\lambda^+:\varphi_\alpha(x,\bar y) = \varphi_*(x,\bar y)\}$
 is a set of cardinality $\lambda^+$, so it is enough to prove 
that for every finite non-empty
$w \subseteq W$ there is $b_w \in
\cap\{\varphi_\alpha(\gC,a_\alpha):\alpha \in W \backslash w\} \cap
\bigcap\{\varphi_*(\gC,\bar a_\alpha):\alpha \in w \backslash
\{\max(w),\min(w)\}$.

Now letting $\alpha_0 < \ldots < \alpha_n$ list $w$, we can choose
$b_w := b_{\alpha_0} \ldots b_{\alpha_n}$ (the product).  Now on the one hand,
$b_w \in \varphi_*(\gC,\bar a_\alpha)$ for $\alpha \in W \backslash w$
because $\varphi_*(\gC,\bar a_\alpha)$ by the choice of $\varphi_\alpha$
 is closed under products and
$\ell < n \Rightarrow b_{\alpha_\ell} \in \varphi_*(\gC,\bar a_\alpha)$.
On the other hand if $\ell = 1,\dotsc,n-1$ and $b_u \in
\varphi_*(\gC,\bar a_{\alpha_\ell})$ by the ``moreover" in clause (b)
of $\circledast$ we are done.
\end{PROOF}
\newpage

\section {$\bbL_{\infty,\bar\kappa}(\tau_T)$-definable
subgroups of bounded index} 

The main result of this section is Theorem \ref{h.24}: if the monster model
${\gC}$ is $\bar\kappa$-saturated, $T = Th({\gC})$ is dependent, $G$
is a definable abelian group over $A_*,|A_*| < \bar\kappa$ \then \,
$H_* = \cap\{H:H$ is a union-type-definable
\footnote{recall that this means: is of the form $\cup\{p_i({\gC}):i <
i_*\}$ where for some $A \subset {\gC}$ of cardinality $<
\bar\kappa$, each $p_i$ is a type over $A$.} subgroup of $G$ of bounded index,
i.e. $< \bar\kappa\}$ \then \, $H_*$ is a union-type-definable over $A_*$ and
has index $\le 2^{|T|+|A_*|}$.  As
usual $\bar\kappa$ is strongly inaccessible (strong limit of large
enough cofinality such that $(\cH(\bar\kappa),\in) \prec_{\Sigma_n} (\mathbf
V,\in)$ some $n$ large enough is enough).
 
\begin{context}
\label{h.0}
\mn
\begin{enumerate}
\item[$(a)$]    ${\gC}$ is a monster ($\bar \kappa$-saturated)
model of the complete first order theory $T$

(we assume $\bar\kappa$ is strongly inaccessible $> |T|$; this just
for convenience, we do not really need to assume there is such cardinal)
\sn
\item[$(b)$]   $p_*(x)$ is a type and $(x,y) \mapsto x * y,x \mapsto
x^{-1}$ and $e_G$ are first order definable (in ${\gC}$)  
two-place function, one place function and element
(with parameters $\subseteq \text{ Dom}(p_*)$ for
simplicity) such that their restriction to $p_*({\gC})$ gives
it a group structure which we denote by $G =
G^{\gC}_{p_*}$.  Let Dom$(p_*) = A_*$; we may write $ab$ instead of
$a * b$.  When the group is Abelian we may use the additive notation.
\end{enumerate}
\end{context}

\begin{notation}
\label{h.1}
For an $m$-type $p(\bar x)$ let 
$p(B) = \{\bar b \in {}^{\ell g(\bar x)}B:\bar b$ realizes 
$p(\bar x)\}$, so if $p = p(x)$ then we stipulate $p({\gC}) \subseteq 
{\gC}$.  For a set
$\mathbf P$ of $< \bar\kappa$ of $m$-types let $\mathbf P(B) =
\cup\{p(B):p \in \mathbf P\}$.
\end{notation}

\begin{definition}
\label{h.4}
For any $\alpha < \bar\kappa$ and sequence $\bar{\mathbf a} 
\in {}^\alpha{\gC}$ and $n < \omega$
we define $q^n_{\bar{\mathbf a}}(x) = q^n_{\Gamma_{\bar{\mathbf a}}}(x)$ where
\mn
\begin{enumerate}
\item[$(a)$]   $\Gamma_{\bar{\mathbf a}} = 
\{\varphi(x,\bar{\mathbf a}):\varphi \in \bbL(\tau_T)\}$ and
\sn
\item[$(b)$]    for $\Gamma_1,\Gamma_2$ sets of formulas (possibly
  with parameters) with one free variable $x$ we let:

\begin{equation*}
\begin{array}{clcr}
q^n_{\Gamma_1,\Gamma_2}(x) =: \bigl\{(\exists y_0,\dotsc,y_{2n-1})
[\bigwedge\limits_{\ell < 2n} \,
\psi(y_\ell) &\wedge \bigwedge\limits_{\substack{\ell < n \\ k< m}} 
(\varphi_k(y_{2 \ell},\bar a) \equiv \varphi_k(y_{2 \ell+1},\bar a)) \\
  &\wedge x = (y^{-1}_0 * y_1) * \dotsc * 
(y^{-1}_{2n-2} * y_{2n-1})]: \\
  &(\alpha) \quad \psi(x) \text{ is a finite conjunction with} \\ 
  &\hskip25pt \text{no repetitions of members of } \Gamma_1 \cup p_*(x) \\
  &(\beta) \quad m < \omega \text{ and } \varphi_k(x,\bar a) \in
\Gamma_2 \text{ is a formula for } k < m \\
  & \qquad (\text{ so if } \Gamma_2 = \Gamma_{\bar{\mathbf a}}, \text{ we may
  replace } \bar a \text{ by } \bar{\mathbf a})\bigr\}
\end{array}
\end{equation*}
\sn
\item[$(c)$]   if $\Gamma_1 \subseteq p_*$ we may omit it.  If
$\Gamma_2 = \Gamma_{\bar{\mathbf a}}$ and $\Gamma_1 \subseteq p_*$ 
we may write $q^n_{\bar{\mathbf a}}$ 
instead of $q^n_{\Gamma_{\bar{\mathbf a}}}$.
\end{enumerate}
\end{definition}

\begin{remark}
We have used $x=w(y^{-1}_0 y_1,y^{-1}_2 y_3,\dotsc,y^{-1}_{2n-2} y_{2n-1})$
with $w(x_0,\dotsc,x_{n-1})$ the ``word" $x_0 \ldots x_{n-1}$; but we
may replace $x_0 x_1 \ldots x_{n-1}$ by any group word
$w(x_0,\dotsc,x_{n-1})$ and then may 
write $q^w_{\Gamma_1,\Gamma_2},q^w_\Gamma(x),q^w_{\bar{\mathbf a}}(x)$ instead of
$q^n_{\Gamma_1,\Gamma_2},q^n_{\Gamma_2},q^n_{\bar{\mathbf a}}(x)$
respectively.
Of course, the subgroup of $G$ generated by $q^1_{\bar{\mathbf a}}({\gC})$
includes $q^w_{\bar{\mathbf a}}({\mathfrak C})$ for any non-trivial group word $w$.  
\end{remark}

\begin{observation}
\label{h.8}
For $n < \omega$ and $\bar{\mathbf a} \in {}^\alpha {\gC}$ as above
\mn
\begin{enumerate}
\item[$(a)$]   $q^n_{\bar{\mathbf a}}(x)$ is a 1-type
\sn
\item[$(b)$]   $q^n_{\bar{\mathbf a}}({\gC}) =
\{c_0 \dotsc c_{n-1}$: there are $d_\ell \in p_*({\gC})$ for
$\ell < 2n$ such that $c_\ell = d^{-1}_{2 \ell} d_{2 \ell+1}$ and
$d_{2 \ell},d_{2 \ell +1}$ realize the same type over $\bar{\mathbf a}$
for every $\ell < n\}$
\sn
\item[$(c)$]   $|q^n_{\bar{\mathbf a}}|$ is equal to $|T| + \aleph_0 +
|\ell g(\bar{\mathbf a})| + |p_*(x)|$
\sn
\item[$(d)$]   $q^n_{\bar{\mathbf a}} = \cup\{q^n_{\bar{\mathbf a}
\restriction u}:u \subseteq \ell g(\bar{\mathbf a})$ finite$\}$
\sn
\item[$(e)$]   $q^n_\Gamma = \cup\{q^n_{\Gamma_1}:\Gamma_1 \subseteq
\Gamma$ finite$\}$.
\end{enumerate}
\end{observation}

\begin{PROOF}{\ref{h.8}}
Straightforward. 
\end{PROOF}

\begin{observation}
\label{h.16} 
1) $e_G \in q^n_\Gamma({\gC}) \subseteq q^{n+1}_\Gamma({\gC})$.

\noindent
2) $q^n_\Gamma({\gC})$ is closed under $(-)^{-1}$ and $e_G$
belongs to it so if $a \in q^k_\Gamma({\gC})$ \then \, $a^{-1} \in
   q^k_\Gamma({\gC})$.

\noindent
3) If $a_\ell \in q^{k(\ell)}_\Gamma({\gC})$ for $\ell=1,2$ \then \,
$a_1 * a_2 \in q^{k(1)+k(2)}_\Gamma({\gC})$. 

\noindent
4) $\cup\{q^n_\Gamma({\gC}):n<\omega\}$ is a subgroup of $G$.

\noindent
5) We have $q^m_{\bar{\mathbf a}}({\gC}) \supseteq 
q^n_{\bar{\mathbf b}}({\gC})$ when $n \le m < \omega$ and
Rang$(\bar{\mathbf a}) \subseteq \text{ Rang}(\bar{\mathbf b})$.  
\end{observation}

\begin{PROOF}{\ref{h.16}}
Easy.  
\end{PROOF}

\begin{observation}
\label{h.12}  
1) If $\bar{\mathbf a}_\ell \in {}^{\alpha(\ell)} {\gC}$ for 
$\ell=1,2$ and Rang$(\bar{\mathbf a}_1)
\subseteq \text{ Rang}(\bar{\mathbf a}_2)$ \then \, 
$q^n_{\bar{\mathbf a}_2} \vdash q^n_{\bar{\mathbf a}_1}$.

\noindent
2) Assume $n < \omega$ and $H = p(\gC)$ 
is a subgroup of $G$ of bounded index (i.e. $(G:H) < \bar\kappa$) and
$\alpha < \bar\kappa,\bar{\mathbf a} \in {}^\alpha{\gC}$
has representatives from each left $H$-coset and Rang$(\bar{\mathbf a})$
include $A_*$.  If $p(x)$ is a 1-type $\supseteq p_*(x)$ which is over 
Rang$(\bar{\mathbf a}) = M \prec \gC$ and $p({\gC}) = H$ \then \,
$q^n_{\bar{\mathbf a}}(x) \vdash p(x)$.

\noindent
3) In part (1) if $\mathbf P$ is a set of $(< \bar\kappa)$ 1-types
   $\supseteq p_*(x)$ which are over Rang$(\bar{\mathbf a}) = M \prec
   \gC$ and $\mathbf P({\gC}) = H$ \then  \, $q^n_{\bar{\mathbf a}}({\gC})
   \subseteq H$.

\noindent
4) If $\bar{\mathbf a} \in {}^\alpha {\gC}$ so $\alpha < \bar\kappa$
\then \, we have $\bigcup\limits_{k < \omega} q^k_{\bar{\mathbf a}} 
({\gC})$ is a union-type definable subgroup of $G$ of 
bounded index, in fact of index $\le 2^{|T|+|\alpha|}$.

\noindent
5) If $H$ is union-type-definable over $\bar{\mathbf a}$ where 
$\bar{\mathbf a} \in {}^\alpha {\gC},
\alpha = \ell g(\bar y) < \bar\kappa$ and $H$ 
is a subgroup of $G$ with bounded index, \then \, $(G:H) \le 2^{|T|+|\alpha|}$.
\end{observation}

\begin{claim}
\label{h.13.21}
Any subgroup of $G$ which is
union-type-definable (so using all together $< \bar\kappa$ parameters)
 and is of bounded index (in $G$, bounded means $< \bar\kappa$) 
contains a subgroup of the form $\bigcup\limits_{n} 
q^n_{\bar{\mathbf a}}({\gC})$ for some $\bar{\mathbf a}
\in {}^{\bar\kappa >}{\gC}$ such that $A_* \subseteq 
\text{\rm Rang}(\bar{\mathbf a})$.
\end{claim}

\begin{PROOF}{\ref{h.13.21}}
Let $\langle a_{2,i}:i < i(*)\rangle$ be a 
set of representatives of the left cosets of $H$ in $G$.  
Let $\bar{\mathbf a}_1$ be such that $H = \bigcup_i p_i({\gC}),p_i(x)$
a type over $\bar{\mathbf a}_1$ and let $\bar{\mathbf a}_2 =
\bar{\mathbf a}_1 \char 94 \langle a_{2,i}:i<i(*)\rangle$.  Clearly if
$a,b \in G$ realize the same type over $\bar{\mathbf a}_2$ then $ab^{-1}
\in H$.

\noindent
By applying \ref{h.16} + \ref{h.12}(1) we are done. 
\end{PROOF}

\begin{mc}
\label{h.20}
Assume $T$ is dependent and $G$ is Abelian.  If 
$k(1)+2 < k(2) < \omega$ and $\alpha < \bar\kappa$,
\then \, there are $\lambda < \bar\kappa$ and 
$\bar{\mathbf a}_\varepsilon \in {}^\alpha {\gC}$ for 
$\varepsilon < \lambda$ such that for every 
$\bar{\mathbf a} \in {}^\alpha{\gC}$ satisfying $A_* \subseteq
 \text{\rm Rang}(\bar{\mathbf a})$ we have $\cap\{q^{k(1)}_{\bar{\mathbf
a}_\varepsilon}({\gC}):\varepsilon < \lambda\} \subseteq
q^{k(2)}_{\bar{\mathbf a}}({\gC})$, i.e. $\bigcup\limits_{\varepsilon < \lambda} 
q^{k(1)}_{\bar{\mathbf a}_\varepsilon} \vdash q^{k(2)}_{\bar{\mathbf a}}$.
\end{mc}

\begin{remark}
\label{h.20.4}
1) We may choose $\lambda$ and $\bar{\mathbf a}_\varepsilon$ 
(for $\varepsilon < \lambda$) independently of $k(1),k(2)$.

\noindent
2) The proof says that $\lambda = 2^{2^{|\alpha|+|T|+|A_*|}}$ is enough.
\end{remark}

\begin{question}  Is $G$ abelian or the subgroup abelian necessary?
\end{question}

\begin{PROOF}{\ref{h.20}}
\underline{Stage A}:   Let $\lambda := 
(2^{2^{|\alpha|+|T|+|A_*|}})$, and assume that the desired conclusion 
fails for $\lambda$.
\medskip

\noindent
\underline{Stage B}: As this fails, we can find $\langle \bar{\mathbf
a}_\varepsilon,c_\varepsilon:\varepsilon < \lambda^+\rangle$ such that:
\mn
\begin{enumerate}
\item[$(a)$]   $\bar{\mathbf a}_\varepsilon \in {}^\alpha {\gC}$ and $A_*
\subseteq \text{ Rang}(\bar a_\varepsilon)$
\sn
\item[$(b)$]   $c_\varepsilon \in 
\bigcap\{q^{k(1)}_{\bar{\mathbf a}_\zeta}({\gC}):
\zeta < \varepsilon\} \subseteq G$ equivalently
$c_\varepsilon \in q^{k(1)}_{\bar{\mathbf a}_\zeta}({\gC})$ for
$\zeta < \varepsilon$
\sn
\item[$(c)$]   $c_\varepsilon \notin 
q^{k(2)}_{\bar{\mathbf a}_\varepsilon}({\gC})$. 
\end{enumerate}
\mn
[Why?  Choose $(\bar{\mathbf a}_\varepsilon,c_\varepsilon)$ by induction
on $\varepsilon < \lambda^+$.]
\medskip

\noindent
\underline{Stage C}:  Without loss of generality $\langle 
(\bar{\mathbf a}_\varepsilon,c_\varepsilon):\varepsilon< \lambda^+\rangle$ is an
indiscernible sequence over $A_*$.

\noindent
[Why?  We first use Erd\"os-Rado theorem to repalce clause (c) by
\mn
\begin{enumerate}
\item[$(c)'$]  $c_\varp \in (\neg \bigwedge\limits_{\ell < m} 
\varphi_\ell)(\gC,\bar{\mathbf a}_\varp)$ 
for some $m$ and $\varphi_\ell(x,\bar{\mathbf a}_\varp) \in
q^{k(2)}_{\bar{\mathbf a}_\varp}$
\end{enumerate}
\mn
and then use Ramsey theorem (and compactness, i.e. saturation of $\gC$.]
\medskip

\noindent
\underline{Stage D}:  For $\varepsilon < \lambda^+$ let 
$c'_\varepsilon = c^{-1}_{2 \varepsilon} \ast c_{2 \varepsilon+1}$ and
let $\bar{\mathbf a}'_\varepsilon = 
\bar{\mathbf a}_{2 \varepsilon}$ and for any finite 
$u \subseteq \lambda^+$ let $c'_u = c'_{\varepsilon_0} *
\ldots * c'_{\varepsilon_{m-1}}$ when $\varepsilon_0 < \varepsilon_1 <
\ldots < \varepsilon_{m-1}$ list $u$.
\medskip

\noindent
\underline{Stage E}:  If $u \subseteq \lambda^+$ is finite and $\varepsilon
\in \lambda^+ \backslash u$ then $c'_u$ realizes 
$q^1_{\bar{\mathbf a}'_\varepsilon}$.

\noindent
[Why?  As $G$ is abelian, using additive notation
$c'_u = \Sigma\{c_{2 \zeta +1}:\zeta
\in u\}-\Sigma\{c_{2 \zeta}:\zeta \in u\}$.  By the indiscernibility
the sequences $\langle c_{2 \zeta +1}:\zeta \in u\rangle,
\langle c_{2 \zeta}:\zeta \in u\rangle$ realize the same type over
$\bar{\mathbf a}'_\varepsilon = \bar{\mathbf a}_{2 \varepsilon}$ as
$\varepsilon \notin u$, (noting the specific sequences we use) hence
$\Sigma\{c_{2 \zeta +1}:\zeta \in u\}$ and $\Sigma\{c_{2 \zeta}:\zeta
\in u\}$ are members of $G$ realizing the same type over $\bar{\mathbf
a}'_\varepsilon = \bar{\mathbf a}_{2 \varepsilon}$ hence $c'_u \in
q^1_{\bar{\mathbf a}_\varepsilon}({\gC})$.  So the conclusion is clear.]
\medskip

\noindent
\underline{Stage F}: There is a finite sequence $\bar \vartheta = \langle
\vartheta_\ell(x,\bar{\mathbf y}):\ell < \ell(*)\rangle$ of formulas, $\ell
g(\bar{\mathbf y}) = \alpha$ and $\bar\varphi = \langle \varphi_i(x,\bar
b):i < i(*)\rangle,\varphi_i(x,\bar b) \in p_*(x)$ for $i<i(*) <
\omega,\bar b \subseteq A_*$ such that: for $\varepsilon < \lambda^+$, for no
$d_0,\dotsc,d_{2k(2)-1} \in \bigcap\limits_{i<i(*)} \varphi_i(\gC,\bar
b) \subseteq G$ do we have $c_\varepsilon = 
\sum\limits_{\ell < k(2)} ((-d_{2 \ell}) + d_{2 \ell +1})$ and 
$\ell < \ell(*) \wedge k < k(2) \Rightarrow {\gC} \models
\vartheta_\ell(d_{2k},\bar{\mathbf a}_\varepsilon) \equiv
\vartheta_\ell(d_{2k+1},\bar{\mathbf a}_\varepsilon)$, in other words 
$c_\varepsilon \notin \vartheta^*(\gC,\bar{\mathbf a}'_\varp)$ where
$\vartheta^*(x,\mathbf a'_\varp) \in q^{k(2)}_{\{\vartheta_\ell
(x,\bar{\mathbf a}_\varepsilon):\ell < \ell(*)\}}({\gC})$ says the
above, i.e. is $\neg(\exists
z_0,\dotsc,z_{2k(c)-1})[(\bigwedge\limits_{k<2k(2)} \,
\bigwedge\limits_{i<i(*)} \varphi_i(z_k,b)) \wedge
\bigwedge\limits_{\ell < k(2)} (\vartheta_\ell(z_{2k},\bar{\mathbf
  a}'_\varp) \equiv \vartheta_\ell(z_{2k+1},\bar{\mathbf a}'_\varp)]$. 

\noindent
[Why?  Because $c_\varepsilon \notin q^{k(2)}_{\bar{\mathbf
a}_\varepsilon}({\gC})$, see (c) of stage (B). 
By compactness, we can choose finitely many formulas; by the
indiscernibility their choice does not depend on $\varepsilon$, see 
stage (C).]
\medskip

\noindent
\underline{Stage G}:  If $u \subseteq \lambda^+$ is finite and $\varepsilon
\in u$ \then \, $c'_u \notin q^{k(*)}_{\{\vartheta_\ell(x,\bar{\mathbf
a}'_\varepsilon):\ell < \ell(*)\}} ({\gC})$; moreover, $\bar c'_u
\notin \vartheta^*(\gC,\bar{\mathbf a}')$ for any $k(*) \le k(2) - 1-k(1)$.

\noindent
[Why?  Because
\mn
\begin{enumerate}
\item[$(a)$]   $c'_u = \Sigma\{c'_\zeta:\zeta \in u \backslash
\{\varepsilon\}\} - c_{2 \varepsilon} + c_{2 \varepsilon +1}$
\end{enumerate}
\mn
hence
\mn
\begin{enumerate}
\item[$(b)$]   $c_{2 \varepsilon}  = \Sigma\{c'_\zeta:\zeta \in u \backslash
\{\varepsilon\}\} + c_{2 \varepsilon +1} - c'_u$
\sn
\item[$(c)$]   $\Sigma\{c'_\zeta:\zeta \in u \backslash
\{\varepsilon\}\} \in q^1_{\{\vartheta_\ell(y,\bar{\mathbf a}'_\varepsilon):\ell
< \ell(*)\}}({\gC})$.
\end{enumerate}
\mn
[Why?  As in stage E.]
\mn
\begin{enumerate}
\item[$(d)$]   $c_{2 \varepsilon +1} \in 
q^{k(1)}_{\bar{\mathbf a}_{2\varepsilon}} ({\gC})
\subseteq q^{k(1)}_{\{\vartheta_\ell(x,\bar{\mathbf a}'_\varepsilon):\ell <
\ell(*)\}} ({\gC})$. 
\end{enumerate}
\mn
[Why?  First, the membership by clause (b) of Stage B.  Second, the
inclusion by monotonicity as $\bar{\mathbf a}'_\varepsilon = 
\bar{\mathbf a}_{2 \varepsilon}$.]
\mn
\begin{enumerate}
\item[$(e)$]   $c_{2 \varepsilon} \notin
q^{k(2)}_{\{\vartheta_\ell(x,\bar{\mathbf a}_{2\varepsilon}):
\ell < \ell(*)\}} ({\gC}) = q^{k(2)}_{\{\vartheta_\ell(x,\bar{\mathbf
a}'_\varepsilon):\ell < \ell(*)\}} ({\gC})$; moreover, $c_{2 \varp}
\notin \vartheta^*(\gC,\bar{\mathbf a}_{2 \varp}) =
\vartheta^*(\gC,\bar{\mathbf a}'_\varp)$.
\end{enumerate}
\mn
[Why?  By Stage F.]

Now consider the right side in clause (b).

The first summand is from 
$q^1_{\{\vartheta_\ell(x,\bar{\mathbf a}'_\varepsilon):\ell < \ell(*)\}}
({\mathfrak C})$ by clause (c).

The second summand is from $q^{k(1)}_{\{\vartheta_\ell(x,\bar
a'_\varepsilon):\ell < \ell(*)\}} ({\gC})$ by clause (d).

\noindent
For awhile assume 
$c'_u \in \vartheta^*(\gC,\bar{\mathbf a}'_\varp)$ so $(-c'_u) \in
q^{k(*)}_{\{\vartheta_\ell(x,\bar{\mathbf a}'_\varepsilon):\ell <
\ell(*)\}}({\gC})$ then by \ref{h.16}(3) and the previous two sentence
we deduce $\Sigma\{c'_\zeta:\zeta \in u \backslash \{\varepsilon\}\} +
c_{2 \varepsilon+1} + (-c'_u) \in
q^{k(1)+1+k(*)}_{\{\vartheta_\ell(x,\bar{\mathbf a}'_\varepsilon):\ell <
\ell(*)\}}({\gC})$ which by clause (b) means 
$c_{2 \varepsilon} \in q^{k(1)+1+k(*)}_{\{\vartheta_\ell
(x,\bar a'_\varepsilon):\ell < \ell(*)\}}({\gC})$.  
Now by the assumption of the stage we have
$k(2) \ge k(1) + 1 +k(*)$ but $q^n_\Gamma({\gC})$ does 
$\subseteq$-increase with $n$ hence from the previous sentence
we get a contradiction to clause $(e)$.  So necessarily $c'_u \notin
\vartheta^*(\gC,\bar{\mathbf a}_\varp)$, the desired conclusion of this stage.]
\medskip

\noindent
\underline{Stage H}:  Let $k(*) := k(2)-1-k(1)$.

By Stages (E)+(G) we have: for finite $u \subseteq \lambda^+,c'_u$
realizes $q^{k(*)}_{\{\vartheta_\ell(x,\bar{\mathbf a}'_\varepsilon):\ell <
\ell(*)\}}$ iff $\varepsilon \in u$ iff ${\gC} \models
\vartheta^*(c'_u,\bar{\mathbf a}'_\varepsilon)$.  So
$\vartheta^*(x,\bar{\mathbf y}) \in \bbL(\tau_T)$ has the independence
property, contradiction.  
\end{PROOF}

\begin{theorem}
\label{h.24}  
Assume $T$ is dependent and $G$ is Abelian.
There is $\mathbf P \subseteq \mathbf S(A_*)$ such that:
\mn
\begin{enumerate}
\item[$(a)$]    $\mathbf P({\gC}) := \cup\{p({\gC}):p \in \mathbf P\}$ 
is a subgroup of $G$ of bounded index, i.e. $(G:\mathbf P({\gC})) < \bar\kappa$
\sn
\item[$(b)$]    $(G:\mathbf P(G)) \le 2^{|A_*|+|T|+\aleph_0}$
\sn
\item[$(c)$]    $\mathbf P({\gC})$ is minimal, i.e.,
if $\alpha < \bar\kappa,\bar{\mathbf a} \in {}^\alpha{\gC},
\psi(x,\bar{\mathbf y}) \in \bbL_{\infty,\kappa}(\tau_T)$ and 
$\psi({\gC},\bar{\mathbf a})$ is a subgroup of $G$ \then \, 
$\psi({\gC},\bar{\mathbf a})$ is of bounded index iff $\psi({\gC},\bar{\mathbf a}) 
\supseteq \mathbf P({\gC})$.
\end{enumerate}
\end{theorem}

\begin{PROOF}{\ref{h.24}}
Note:
\mn
\begin{enumerate}
\item[$\circledast_1$]   there is $\mathbf P \subseteq \mathbf S(A_*)$
such that $\mathbf P({\gC})$ is a subgroup of $G$ of bounded index.
\end{enumerate}
\mn
[Why?  Use $\mathbf P = \{p \in \mathbf S(A_*):p_* \subseteq p\}$.]
\mn
\begin{enumerate}
\item[$\circledast_2$]    the family of $\mathbf P$'s as in
$\circledast_1$ is closed under intersection.
\end{enumerate}
\mn
[Why?  As $\bar\kappa$ is strongly inaccessible or just $\cf(\bar\kappa) >
\beth_3(|A_*| +|T| + \aleph_1) \ge 2^{|\{\mathbf P:\mathbf P \subseteq
\mathbf S(A_*)\}|}$ hence the product of $\le 2^{2^{|A_*|+|T|+\aleph_0}}$
cardinals $< \bar \kappa$ is $<\bar\kappa$.]
\mn
\begin{enumerate} 
\item[$\circledast_3$]    Let $\mathbf P^* = \cap\{\mathbf P \subseteq
\mathbf S(A_*):\mathbf P({\gC})$ is a subgroup of a bounded index$\}$.
\end{enumerate}
\mn
We shall show that $\mathbf P^*$ is as required.

So by $\circledast_2 + \circledast_3$
\mn
\begin{enumerate}
\item[$\circledast_4$]   clause (a) holds, i.e.  
$\mathbf P^* \subseteq \mathbf S(A_*)$ and $\mathbf P^*({\gC})$ is a 
subgroup of $G$ of bounded index
\sn
\item[$\circledast_5$]   clause (b), i.e. $(G:\mathbf P^*({\gC})) \le
2^{|A_*|+|T|+\aleph_0}$.
\end{enumerate}
\mn
[Why?  Follows from \ref{h.12}(5).]

Recall
\mn
\begin{enumerate}
\item[$\circledast_6$]   \underline{if} 
$1 \le k < \omega$ and $\mathbf P$ is a set one-types 
(of cardinality $< \bar\kappa$) and
$\mathbf P({\gC})$ a subgroup of $G$ of bounded index 
\then \, $q^k_{\bar{\mathbf a}}({\gC}) \subseteq \mathbf P({\gC})$ 
for some $\bar{\mathbf a} \in {}^\alpha \gC$ for some $\alpha < \kappa$.
\end{enumerate}
\mn
[Why?  See above \ref{h.13.21}.]

Fix $\alpha < \bar \kappa$ and we shall prove that:
\mn
\begin{enumerate}
\item[$\boxtimes$]    if $\bar{\mathbf a}^* \in {}^\alpha{\gC},A_*
\subseteq \Rang(\bar{\mathbf a}_*)$ and $\mathbf P_{\bar{\mathbf a}^*} 
:= \{p \in \mathbf S(\bar{\mathbf a}^*):p$ extend
$q^n_{\bar{\mathbf a}^*}$ for some $n\}$, so $\mathbf P_{\mathbf a^*}({\gC}) = 
\bigcup\limits_{n < \omega} q^n_{\bar{\mathbf a}^*}({\gC})$ is 
a subgroup of $G$ of bounded index \then \, ${\mathbf P}^*({\gC})
\subseteq \mathbf P_{\mathbf a^*}({\gC})$.
\end{enumerate}
\mn
This clearly suffices by $\circledast_4, \circledast_5$ and
$\circledast_6$, i.e.
$\boxtimes$ means that clause (c) of the conclusion holds by
$\circledast_6$.  

Now comes the real point: 
\mn
\begin{enumerate}
\item[$\circledast_7$]    for $k < \omega$ we can choose $\lambda_k
< \bar\kappa$ and $\langle \bar{\mathbf a}^k_\varepsilon:
\varepsilon < \lambda_k \rangle$ such that:
\sn
\item[${{}}$]   $(a) \quad \bar{\mathbf a}^k_\varepsilon \in
{}^\alpha {\gC}$ for $\varepsilon < \lambda_k$
\sn
\item[${{}}$]  $(b) \quad$ for every $\bar{\mathbf a} \in
{}^\alpha {\gC}$ we have $(\bigcup\limits_{\varepsilon < \lambda_k} 
q^k_{\bar{\mathbf a}^k_\varepsilon}) \vdash q^{k+3}_{\bar{\mathbf a}}$.
\end{enumerate}
\mn
[Why is there such a sequence?  By the main claim \ref{h.20} so
actually $\lambda_k = (2^{2^{|T|+|A_*|+|\ell g(\bar a^*)|}})$ suffice
by the proof of \ref{h.20}.]

Define 
\mn
\begin{enumerate}
\item[$\odot_1$]  $X_k := \cap\{q^k_{\bar{\mathbf a}^\varepsilon_k}
({\gC}):\varepsilon < \lambda_k\}$
\sn
\item[$\odot_2$]   $Y_k := \cap\{q^k_{\bar{\mathbf a}}({\gC}):
\bar{\mathbf a} \in {}^\alpha({\gC})\}$.
\end{enumerate}
\mn
Then:
\mn
\begin{enumerate}
\item[$\circledast_8$]   $Y_k \subseteq X_k$.
\end{enumerate}
\mn
[Why?  As $\mathbf a^\varepsilon_k \in {}^\alpha {\gC}$.]
\mn
\begin{enumerate}
\item[$\circledast_9$]   $X_k \subseteq Y_{k+3}$.
\end{enumerate}
\mn
[Why?  By $\circledast_7$.]
\mn
\begin{enumerate}
\item[$\circledast_{10}$]   $\bigcup\limits_{k < \omega} X_k = 
\bigcup\limits_{k < \omega} Y_k$.
\end{enumerate}
\mn
[Why?  By $\circledast_8 + \circledast_9$.]
\mn
\begin{enumerate}
\item[$\circledast_{11}$]   $Y_k = \mathbf P_k({\gC})$ for some
$\mathbf P_k \subseteq \mathbf S(A_*)$.
\end{enumerate}
\mn
[Why?  As any automorphism $F$ of ${\gC}$ over $A_*$ maps $Y_k$
onto itself as it maps $q^k_{\bar{\mathbf a}}({\gC})$ to
$q^k_{F(\bar{\mathbf a})}({\gC})$.]
\mn
\begin{enumerate}
\item[$\circledast_{12}$]   $\bigcup\limits_{k < \omega} Y_k$ is 
$(\bigcup\limits_{k < \omega} \mathbf P_k)({\gC})$.
\end{enumerate}
\mn
[Why?  By $\circledast_{11}$.]
\mn
\begin{enumerate}
\item[$\circledast_{13}$]   $\bigcup\limits_{k < \omega} X_k$ is a subgroup
of $G$ of bounded index.
\end{enumerate}
\mn
[Why?  By \ref{h.16}(3) + \ref{h.12}(4).]

Recall $\mathbf P_{\bar{\mathbf a}_*}$ is from $\boxtimes$ above.
\mn
\begin{enumerate}
\item[$\circledast_{14}$]   $Y_k \subseteq \mathbf P_{\bar{\mathbf a}^*}(\gC)$.
\end{enumerate}
\mn
[Why?  By the definition of $Y_k$ and $\circledast_7$ we have $Y_k
\subseteq q^k_{\bar{\mathbf a}^*}({\gC}) \subseteq 
\mathbf P_{\bar{\mathbf a}^*}({\gC})$.]

Let us sum up and prove $\boxtimes$ thus finishing: 
$\mathbf P_{\bar{\mathbf a}^*}(\gC)$ include $\bigcup\limits_{k < \omega} Y_k$ by
$\circledast_{14}$ and $\bigcup\limits_{k < \omega} Y_k$ is equal to 
$\bigcup\limits_{k < \omega} X_k$ by $\circledast_{10}$, and is 
equal to $(\bigcup\limits_{k < \omega} \mathbf P_k)({\gC})$ 
by $\circledast_{12}$.
Hence by $\circledast_{13}$ we know that $(\bigcup\limits_{k < \omega} \mathbf
P_k)({\gC})$ is a subgroup of $G$ of bounded index, hence by the
definition of $\mathbf P^*$ in $\circledast_3$ we know that $\bigcup\limits_{k <
\omega} \mathbf P_k \subseteq \mathbf P^*$.  Hence $(\bigcup\limits_{k < \omega}
\mathbf P_k)({\gC})$ includes $\mathbf P^*({\gC})$.
So $\mathbf P_{\bar{\mathbf a}^*}({\gC}) \supseteq 
\bigcup\limits_{k < \omega} Y_k =
(\bigcup\limits_{k < \omega} \mathbf P_k)({\gC}) 
\supseteq \mathbf P^*({\gC})$ as required in $\boxtimes$.

So we have proved $\boxtimes$ hence has proved the conclusion.
\end{PROOF}
\newpage

\section {On 2-dependent $T$}

We try to see what, from $2$-dependence of $T$, we can deduce on
definable groups.  On $n$-dependent $T$ see \cite[\S5 (H)]{Sh:863}.

\begin{hypothesis}
\label{dt.7}
\mn
\begin{enumerate}
\item[$(a)$]    $T$ be first order complete, ${\gC} = {\gC}_T$
\sn
\item[$(b)$]   $G$ is a type definable group over $A_*$, i.e. for some
  1-type $p_*,G$ has a set of elements
$p_*({\gC})$ and the functions $(x,y) \mapsto xy,x \mapsto x^{-1},e_G$
which are definable over $A_* = \text{ Dom}(p_*)$; this is irrelevant for
\ref{dt.22} - \ref{dt.25}.
\end{enumerate}
\end{hypothesis}

\begin{definition}
\label{dt.14}
For a set $B \subseteq {\gC}$ let
\mn
\begin{enumerate}
\item[$(a)$]   $\mathbf R_B = \{q:q=q(x)$ is a 1-type over $B$ and 
$G_q$ a subgroup of $G$ of index $< \bar \kappa\}$ where 
\sn
\item[$(b)$]   $G_q = G[q] = 
\{a \in G:a$ realizes $q\} = (p_* \cup q)({\gC})$
\sn
\item[$(c)$]   $q_B = q[B] = \cup\{q:q \in \mathbf R_B\}$ and $G_B =
q_B({\gC}) \cap G$  and, of course
\sn
\item[$(d)$]   $\mathbf R_{\bar b} = \mathbf R_{\text{Rang}(\bar b)}$
and $q_{\bar b} = q_{\text{Rang}(\bar b)}$ and $G_{\bar b} =
G_{\text{Rang}(\bar b)}$.
\end{enumerate}
\end{definition}

\begin{observation}
\label{dt.18}
1) $q_B = \cup\{q_{\bar b}:\bar b
\in {}^{\omega \ge}B\} = \cup\{q \in \mathbf R_B:q$ countable$\}$.

\noindent
2) $q_B$ is $\subseteq$-maximal in $\mathbf R_B$.

\noindent
3) $G_B = \cap\{G_q:q \in \mathbf R_B\} = \cap\{G_q:q \in \mathbf R_B$ is
countable$\}$ and is $\subseteq$-minimal in $\{G_q:q \in \mathbf R_B\}$.

\noindent
4) If $q \in \mathbf R_B$ and $q' \subseteq q$ is countable then we can
find a sequence $\langle \psi_n(x,\bar a):n < \omega\rangle$ of finite
conjunctions of members of $q$ such that:
\mn
\begin{enumerate}
\item[$(\alpha)$]  $p_*(x) \cup \{\psi_n(x,\bar a):n < \omega\}
\vdash q'(x)$
\sn
\item[$(\beta)$]   $\psi_{n+1}(x,\bar a) \vdash \psi_n(x,\bar a)$
\sn
\item[$(\gamma)$]   $\bar a \in {}^{\omega \ge} B$, (for notational
simplicity we allow it to be infinite)
\sn
\item[$(\delta)$]   $p_*(x) \cup p_*(y) \cup \{\psi_{n+1}(x,\bar
a),\psi_{n+1}(y,\bar a)\} \vdash \psi_n(x y,\bar a) \wedge
\psi_n(x^{-1},\bar a) \wedge \psi_n(x y^{-1},\bar a)$. 
\end{enumerate}
\mn
5) In part (4), if we allow $\psi_n(x,\bar a_n)$ to be a finite
conjunction of members of $q \cup p_*$ (e.g. if $p_* \subseteq q$)
\then \, we can omit $p_*$ in clauses $(\alpha),(\delta)$ so
$\bigcap\limits_{n < \omega} \psi_n({\gC},\bar a)$ is a group.

\noindent
6) There is a countable $p'(x) \subseteq p_*(x)$ such that 
$p'({\gC})$ is a group under the definable functions $(x,y) \mapsto xy,x
\mapsto x^{-1},e_G$, moreover there is a sequence $\langle
\psi_n(x,\bar a):n < \omega\rangle$ of finite conjunctions of members
of $p_*(x)$ such that:
\mn
\begin{enumerate}
\item[$(\alpha)$]   $\bigwedge\limits_{\ell \le 2} \psi_0(x_\ell,\bar a)
\vdash (x_0 x_1)x_2 = x_0(x_1 x_2) \wedge x_0 e_G = e_G x_0 = x_0 \wedge
x_0 x^{-1}_0 = x^{-1}_0 x_0 = e_G$, (implicitly this means that $x_0
x_1,(x_0 x_1) x_2,x_1 x_2,x_0(x_1,x_2)$ and $x^{-1}_0$ are well
defined)
\sn
\item[$(\beta)$]   $\psi_{n+1}(x,\bar a) \vdash \psi_n(x)$
\sn
\item[$(\gamma)$]  $\psi_{n+1}(x,\bar a) \wedge \psi_{n+1}(y,\bar
a) \vdash \psi_n(xy,\bar a) \wedge \psi_n(x^{-1},\bar a) \wedge
\psi_n(xy^{-1},\bar a)$.
\end{enumerate}
\end{observation}

\begin{PROOF}{\ref{dt.18}}
Obvious and as in \cite{Sh:876}.
\end{PROOF}

\begin{definition}
\label{dt.22}
1) We say $T$ is $2$-independent \when \, we can 
find an independent sequence of formulas of the form
$\langle \varphi(\bar x,\bar b_n,\bar c_m):n,
m < \omega \rangle$ in ${\gC} = {\gC}_T$ or just in some model of $T$.

\noindent
2) $T$ is $1/2$-dependent (or ``$T$ is $2$-dependent") 
means the negation of $2$-independent (see \cite[\S5 (H)]{Sh:863}). 

\noindent
3) We say $\varphi(\bar x,\bar y_0,\dotsc,\bar y_{n-1})$ is
$n$-independent (for $T$) \when \, in ${\gC}_T$ we can, for each
$\lambda < \bar \kappa$, find $\bar a^\ell_\alpha \in {}^{\ell g(\bar
y_\ell)}({\gC}_T)$ for $\alpha < \lambda,\ell < n$ such that the
sequence $\langle \varphi(\bar x,\bar a^0_{\eta(0)},\dotsc,\bar
a^{n-1}_{\eta(n-1)}):\eta \in {}^n \lambda\rangle$ is an independent
sequence of formulas.

\noindent
4) $T$ is $n$-independent \when \, some formula $\varphi(\bar x,\bar
   y_0,\bar y_1,\dotsc,\bar y_{n-1})$ is $n$-independent.

\noindent
5) $T$ is $(1/n)$-dependent (or $T$ is $n$-dependent)
\when \, it is not $n$-independent.
\end{definition}

\begin{remark}
1) In fact $T$ is $n$-independent iff some
$\varphi(x,\bar y_0,\dotsc,\bar y_{n-1})$ is $n$-independent (for
$T$).  We shall write it down in \ref{dt.56} below.

\noindent
2) So 1-independent means independent.
\end{remark}

\begin{claim}
\label{dt.56}
1) For a complete first order theory $T$, there is a 2-independent formula 
$\varphi(x,\bar y,\bar z)$ in $T$ \underline{iff} $T$ is 2-dependent, i.e.
for some $m \ge 1$ there is a 2-independent formula $\varphi(\bar x_m,\bar
y,\bar z)$ with $\bar x_m = \langle x_\ell:\ell <m\rangle$ \underline{iff}
this holds for every $m \ge 1$.

\noindent
2) Similarly for $k$-independent.

\noindent
3) Moreover, if the formula $\varphi(\bar x_m,\bar y_0,\dotsc,
\bar y_{k-1})$ is $k$-independent, \then \, for some $n < m$ and
$\langle b_i:i < m,i \ne m\rangle$ the formula $\varphi =
\varphi(b_0,\dotsc,b_{n-1},x,b_{n+1},\dotsc,b_{m-1},
\bar y_0,\dotsc,\bar y_{k-1})$ is $k$-independent.
\end{claim}

\begin{PROOF}{\ref{dt.56}}  
1) By (2).

\noindent
2) Easily the third statement implies the second, obviously the first
statement implies the third as trivially as we can add dummy
variables.  For the ``second implies the first" direction we 
prove it by induction on
$m$; so assume $k < \omega,\bar x_m = \langle x_\ell:\ell < m\rangle$ and the
formula $\varphi(\bar x_m,\bar y_0,\dotsc,\bar
y_{k-1})$ is $k$-independent.  Let $n_\ell = \ell g(\bar y_\ell),n = 
\sum\limits^{k}_{\ell=1} n_\ell$, of course, \wilog \, $m>1$.  
This means that in
${\gC}_T$ we have $\bar a_{\ell,i} \in 
{}^{\ell g(\bar y_\ell)}{\gC}$ for $\ell <k,i < \omega$ 
such that the sequence $\langle
\varphi(\bar x_m,\bar a_{0,\eta(0)},\dotsc,\bar
a_{k-1,\eta(k-1)}):\eta \in {}^k \omega\rangle$ of formulas 
is independent.  Let
inc$_{n}(\omega) = \{\eta \in {}^n \omega:\eta$ increasing$\}$,
similarly inc$_{<n}(\omega)$.

So for any $R \subseteq \text{ inc}_k(\omega)$ there is $\bar b_R \in
{}^m{\gC}$ such that ${\gC} \models \varphi[\bar b_R,\bar
a_{0,\eta(0)},\dotsc,\bar a_{k-1,\eta(k-1)}]^{\text{if}(\eta \in R)}$
for $\eta \in \text{ inc}_n(\omega)$.

As we can add dummy variables \wilog \, $\bar a_{\ell,i} = \bar a_i$,
i.e. it does not depend on $\ell$ and also $\langle \bar a_i:i <
\omega\rangle$ is an indiscernible sequence.

Let $R_* \subseteq \text{ inc}_k(\omega)$ be random
enough\footnote{which means that e.g. choose a countable $N \prec
  (\cH(\aleph_1),\in)$ and $R_*$ does not belong to any null Borel subset
  of $\cP(\inc_k(\omega))$; the probability space is defined by considering
  $R_*$ as the random subset and the events ``$\eta \in R_*$" for
  $\eta \in \inc_k(\omega)$ are independent, each of probability $1/2$.}.

For $\eta \in \text{ inc}_{< \omega}(\omega)$, let 
$\bar a_\eta = \bar a_{\eta(0)} \char 94 \bar a_{\eta(1)} \char 94
\ldots \char 94 \bar a_{\eta(k-1)}$.  
We say $u_1,u_2 \subseteq \omega$ are $R_*$-similar if
$|u_1| = |u_2|$ and the one-to-one order preserving function $h$ from
$u_1$ onto $u_1$ is an isomorphism from $(u_1,R_* \restriction u_1)$ onto
$(u_2,R_* \restriction u_2)$.  Let $I$ be the model $(\omega,<,R_*)$.

Without loss of generality (by Nessetril-Rodl theorem see
e.g. \cite{GrRoSp}; on such uses see \cite[Ch.III,\S1]{Sh:e} =
\cite{Sh:E59}, in particular \cite[1.26]{Sh:E59})
\mn
\begin{enumerate}
\item[$(*)_1$]   in $\gC$, the sequence $\langle \bar a_t:t \in I\rangle$ is
indiscernible above $\bar b_{R_*}$ which means
\sn
\item[$(*)_2$]   if $j<\omega$ and $\eta_\ell \in 
\text{ inc}_j(\omega),u_\ell = \text{ Rang}(\eta_\ell)$ for $\ell=1,2$ and
$u_1,u_2$ are $R_*$-similar then $\bar a_{\eta_1},\bar
a_{\eta_2}$ realizes the same type in ${\gC}$ over $\bar b_{R_*}$.
\end{enumerate}
\mn
Let $\bar b_{R_*} = \bar b_1 \char 94 \bar b_2$ where $\ell g(\bar
b_1) < m,\ell g(\bar b_2) < m$.

So by $(*)_1$ we have
\mn
\begin{enumerate}
\item[$(*)_3$]   if $\eta \in \text{ inc}_k(\omega)$ \then \, the 
value tp$(\bar a_\eta,\bar b_2,{\mathfrak C})$ depends just on truth value
of $\eta \in R_*$.
\end{enumerate}
\mn
First, assume $\{\text{tp}(\bar a_\eta,\bar b_2,{\gC}):\eta \in 
\inc_k(\omega)\}$ is constant, then letting $\bar x'
= \bar x \rest \ell g(\bar b_1)$ the formula $\varphi(\bar x',
\bar b_2,\bar y_0,\dotsc,\bar y_{k-1})$ is $k$-independent when
$\langle \bar a_i:i < \omega\rangle$ is indiscernible over $\bar b_2$;
pedantically the formula $\varphi(\bar x',(\bar x \rest [\ell g(\bar
b_2) \char 94 \bar y_0),\bar y_1,\dotsc,\bar y_{k-1})$ is
(the restriction to increasing sequences is not serious, 
see \ref{nd.17}(4) below), so by the induction hypothesis we are done.

If $\langle \bar a_i:i < \omega\rangle$ is not an indiscernible
sequence over $\bar b_\eta$, then let $n$ be minimal such that we can choose
$i_0 < \ldots i_{n-1} < \omega,j_0 < \ldots < j_{n-1} < \omega$ we
have $\tp(\bar a_{i_0} \char 94 \ldots \char 94 \bar a_{i_{n-1}},\bar
b_2) \ne \tp(\bar a_{j_0} \char 94 \ldots \char 94 \bar a_{j_{n-1}},\bar
b_2)$.  By $(*)_2$ the sets $\{i_\ell:\ell < n\},\{j_\ell:\ell < n\}$
are not $R_*$-similar.  \Wilog \, there is a unique $v \in [n]^k$ such
that $\langle i_\ell:\ell \in v\rangle \in R_* \leftrightarrow \langle
j_\ell:\ell \in v\rangle \notin R_*$.  Now playing with indiscernible
and using $\bar b_2,\{\bar a_i:i < n,i \in v\}$ as parameters we can finish.

Second assume that $\{\text{tp}(\bar a_\eta,\bar b_2,{\gC}):\eta
\in \text{ inc}_k(\omega)\}$ is not constant, by $(*)_1$ equivalently
$(*)_2$ this means that there is a formula 
$\psi(\bar x_{\ell g(\bar b_1)},\bar y_0,\dotsc,\bar y_{k-1})$ such that:
\mn
\begin{enumerate}
\item[$(*)_4$]   if $\eta \in \text{ inc}_k(\omega)$ then ${\gC}
\models \psi[\bar b_2,\bar a_{\eta(0)},\dotsc,
\bar a_{\eta(k-1)}]$ iff $\eta \in R_*$.
\end{enumerate}
\mn
This also suffices by $\langle \bar a_i:i < \omega\rangle$ being an
indiscernible sequence.

\noindent
3) As in part (2).  
\end{PROOF}

\begin{observation}
\label{nd.17}
1) $T$ is $k$-dependent \when \, : for
every $m,\ell$ and $\varphi(\bar x_m,\bar y) \in 
\bbL(\tau_T)$ for infinitely many $n < \omega$ we have 
$|A| \le n \Rightarrow |\mathbf S^m_{\{\varphi(\bar x_m,\bar y)\}}(A)| 
< 2^{(n/\ell)^k}$.

\noindent
2) In fact we can restrict ourselves to $m=1$ and/or we can replace
$\{\varphi(\bar x_m,\bar y)\}$ by $\Delta = \{\varphi_\ell(\bar
x_m,\bar y_\ell):\ell < \ell_*\}$.

\noindent
3) For any $k,T$ is $k$-independent iff $T^{\text{eq}}$ is $k$-independent.

\noindent
3A) $T$ is not $k$-independent \when \, for every
$\varphi(x,\bar y_0,\dotsc,\bar y_{k-1}) \in \bbL(\tau({\gC}^{\text{eq}}_T))$,
for infinitely many $n$ we have $A \subseteq
{\gC}^{\text{eq}} \wedge |A| \le kn \Rightarrow 
|\mathbf S_{\{\varphi(x,\bar y_0,\dotsc,\bar y_{k-1}\}}(A)| < 2^{n^k}$.

\noindent
4) In Definition \ref{dt.22}(3) we can restrict ourselves to
``increasing $\eta$", similarly in \ref{dt.22}(1).  In fact, 
$\varphi(\bar x,\bar y_0,\dotsc,\bar y_{k-1})$ is $k$-independent
\underline{if} for every $n$ there are $\bar a_{\ell,m} \in {}^{\ell g(\bar
y_\ell)}{\gC}$ for $m < n,\ell < k$ such that $\langle
\varphi(\bar x,\bar a_{0,\eta(0)},\dotsc,
\bar a_{k-1,\eta(k-1)}):\eta \in \inc_k(n)\rangle$ is an
independent sequence of formulas. 
\end{observation}

\begin{PROOF}{\ref{nd.17}}
1) Straightforward and see \cite[\S5 (G)]{Sh:863}.

\noindent
2) Similarly using \ref{dt.56} above.

\noindent
3) Easy by the definition.

\noindent
3A) By parts (1),(2).

\noindent
4) It is enough to prove the second sentence; for every $n$ we
   first find $\langle \bar a_{\ell,m}:m < nk,\ell < k\rangle$ as
   guaranteed there (for $nk$).  Now let 
$\bar a'_{\ell,m} = \bar a_{\ell,\ell n+m}$ so 
$\langle \bar a'_{\ell,m}:m<n,\ell <k\rangle$
   are as required in \ref{dt.22}(3) for $\lambda = n$.  By compactness
   equivalently, by ``${\gC}$ is $\bar\kappa$-saturated" we are done.
\end{PROOF}

\begin{example}
\label{dt.23}
Let $k \ge 1$, a natural $k$-independent but $1/(k+1)$-dependent theory, as
simple as possible, is the model completion of the following theory
(so for $k=1$ this is a (1/2)-dependent, independent $T$):
\mn
\begin{enumerate}
\item[$(A)$]   the vocabulary is

$P_\ell \, (\ell < k+1)$, unary predicates

$R$, a $(k+1)$-place predicate
\sn
\item[$(B)$]   $M$ a model of $T$ \underline{iff}:
\begin{enumerate}
\item[$(a)$]    $\langle P^M_\ell:\ell \le k \rangle$ is a partition
of $|M|$
\sn
\item[$(b)$]    $R^M \subseteq \{(a_0,\dotsc,a_k):a_\ell \in P^M_\ell$
for $\ell=0,\dotsc,k\}$.
\end{enumerate}
\end{enumerate}
\mn
[Note first that clearly the model completion exists and has
elimination of quantifiers.  Second, the formula
$R(x,y_0,\dotsc,y_{k-1})$ exemplifies that $T_k$ is $k$-independent.
Third, $T$ is $1/(k+1)$-dependent by \ref{nd.17}(1) and the
elimination of quantifiers.]
\end{example}

\begin{example}
\label{dt.23c}
Let $T^*_n$ be the theory with the
vocabulary $\{R\},R$ is $n$-place predicate such that $T$ is 
saying $R$ is symmetric and
irreflexive (i.e. $\bigvee\limits_{i<j<n} x_i = x_j
\rightarrow \neg R(x_0,\dots,x_{n-1})$ and $\bigwedge\{R(x_{\varkappa(0)},
\dotsc,x_{\varkappa(n-1)}) \equiv R(x_{\pi(0)},
\dotsc,x_{\pi(n-1)})$ for permutations $\varkappa,\pi$ of $\{0,\dotsc,n-1\}$.

Let $T_n$ be the model completion of $n$.  Then $T_n$ is
$(n-1)$-independent but not $n$-independent (for not
$n$-independent use \ref{nd.17}(1)).
\end{example}

\begin{example}
\label{dt.25}
Any theory $T$ of an infinite Boolean algebras is
(independent, moreover) $k$-independent for every $k$.

\noindent
[Why?  Let for simplicity $k=2$.
Let $\varphi(x,y,z) := (y \cap z \le x)$.  Now for any $n$, let
$\bbB$ be a Boolean sub-algebra of ${\gC}_T$ 
with pairwise disjoint non-zero elements
$\langle a_{i,j}:i,j<n\rangle$ and let $b_i := \cup\{a_{i,j}:j<n\}$
and $c_j := \cup\{a_{i,j}:i<n\}$.  Now $\langle
\varphi(x,b_i,c_j):i,j<n\rangle$ are independent because for $u
\subseteq n \times n$ the element $a_u := \cup\{a_{i,j}:(i,j) \in u\}$
realizes the type $\{\varphi(x,b_i,c_j)^{\text{if}((i,j) \in u)}:
(i,j) \in u\}$.  For $n > 2$ we use $\varphi(x,y_0,\dotsc,y_{n-1}) =
``(y_0 \cap y_1 \cap \ldots \cap y_{n-1}) \le x"$.]
\end{example}

\noindent
Now comes the property concerning the definable group $G$ which interests us.
\begin{definition}
\label{dt.28}
We say that $G$ has
$\kappa$-based bounded subgroups \when \,: for every $\kappa$-saturated
$M \prec \gC$ which include $A_*$ (hence $p_* \subseteq q_M$) and 
$\bar b \in {}^{\omega >}({\gC})$ there is $r \in 
\mathbf R_{M \cup \bar b}$ of cardinality $< \kappa$ such that 
$q_{M \cup \bar b}$ is equivalent to $q_M \cup r$ (equivalently $q_M
\cup r \vdash q_{M \cup \bar b}$), see Definition \ref{dt.14}.
\end{definition}

\noindent
The main result here is
\begin{theorem}
\label{dt.35}  
If $T$ is $1/2$-dependent and $\kappa =
\beth_2(|T| + |p_*|)^+$ or just $\kappa = 
((||T|+|p_*|)^{\aleph_0})^+ + \beth^+_2$
(and $G$ is as in \ref{dt.7}) \then \, $G$ has 
$\kappa$-based bounded subgroups.
\end{theorem}

\begin{PROOF}{\ref{dt.35}}
Assume not and let $\theta = \aleph_0$.  Let $M$ and $\bar b \in
{}^{\omega >}{\gC}$ form a counter-example.  Then we choose the triple
$(r_\alpha,\bar c_\alpha,d_\alpha)$ by induction on $\alpha < \kappa$
such that:
\mn
\begin{enumerate}
\item[$\circledast_1$]   $(a) \quad \bar c_\alpha \in {}^{\omega \ge} M$
\sn
\item[${{}}$]   $(b) \quad r_\alpha = r_\alpha(x,\bar c_\alpha,\bar b) 
=\{\psi^\alpha_n(x,\bar c_\alpha,\bar b):n < \omega\} \in
\mathbf R_{\bar c_\alpha \char 94 \bar b}$, see \ref{dt.14}
\sn
\item[${{}}$]   $(c) \quad d_\alpha \in 
G_{\cup\{\bar c_\beta:\beta < \alpha\} \cup \bar b}$ or just (which
follows) $d_\alpha \in G_{r_\beta}$ for $\beta < \alpha$
\sn
\item[${{}}$]   $(d) \quad d_\alpha \notin G_{r_\alpha}$ moreover
  \wilog \, 
${\mathfrak C} \models \neg \psi^\alpha_0(d_\alpha,\bar c_\alpha,\bar b)$.
\sn
\item[${{}}$]   $(e) \quad d_\alpha \in q_M(\gC)$
\sn
\item[${{}}$]   $(f) \quad \{\psi^\alpha_{n+1}(x,\bar c_\alpha,\bar b)\},
\psi^\alpha_{n+1}(y,\bar c_\alpha,\bar b)\}$

\hskip40pt $\vdash \psi^\alpha_n(xy,\bar c_\alpha,\bar b) 
\wedge \psi^\alpha_n(x^{-1},\bar c_\alpha,\bar b) 
\wedge \psi^\alpha_n(xy^{-1},\bar c_\alpha,\bar b)$.
\end{enumerate}
\mn
[Why we can?  By the assumption toward contradiction.]

Now as cf$(\kappa) > |T|^{\aleph_0}$ \wilog
\mn
\begin{enumerate}
\item[$\circledast_2$]   $\psi^\alpha_n = \psi_n$ for $\alpha < \kappa$.
\end{enumerate}

Of course
\mn
\begin{enumerate}
\item[$\circledast_3$]   $(G:G_{r_\alpha}) \le 2^{\aleph_0}$.
\end{enumerate}
\mn
[Why?  Otherwise let $a_\varepsilon \in G$ for $\varepsilon <
(2^{\aleph_0})^+$ be such that $\langle a_\varepsilon 
G_{r_\alpha}:\varepsilon < (2^{\aleph_0})^+\rangle$ is
without repetition.  For each $\varepsilon < \zeta <
(2^{\aleph_0})^+$ let $n_{\varepsilon,\zeta}$ be the minimal $n$ such
that ${\mathfrak C} \models \neg \psi_n(a^{-1}_\varepsilon a_\zeta,\bar
c_\alpha,\bar b)$, so by Erd\"os-Rado theorem for some 
$n(*)$ and infinite ${\cU} \subseteq
(2^{\aleph_0})^+$ we have $n_{\varepsilon,\zeta} = n(*)$ for
$\varepsilon < \zeta$ from ${\cU}$.  By compactness we can find
$a_\varepsilon \in G$ for $\varepsilon < \bar\kappa$ such that
$\varepsilon < \zeta < \bar \kappa \Rightarrow {\gC} \models \neg
\psi_{n(*)}[a^{-1}_\varepsilon a_\zeta,\bar c_\alpha,\bar b]$,
contradiction to $(G:G_{r_\alpha}) < \bar\kappa$.]
\mn
\begin{enumerate}
\item[$\circledast_4$]   there is ${\cU} \in [\kappa]^\theta$ 
such that: if $\alpha < \beta < \gamma$ are 
from ${\cU}$ then $d^{-1}_\alpha d_\beta \in G_{r_\gamma}$.
\end{enumerate}
\mn
[Why?  For each $\alpha < \kappa$ let $\langle
a_{\alpha,\varepsilon} G_{r_\alpha}:\varepsilon < \varepsilon_\alpha
\le 2^{\aleph_0}\rangle$ be a partition of $G$.  For $\alpha <\beta <
\kappa$ let $\varepsilon = \varepsilon_{\alpha,\beta}$ be such that
$d_\alpha \in a_{\beta,\varepsilon} G_{r_\beta}$.  As $\kappa \rightarrow
(\theta)^2_{2^{\aleph_0}}$ because $\beth^+_2 \rightarrow
(\omega)^2_{2^{\aleph_0}}$ clearly, for some ${\cU} \in [\kappa]^\theta$
and $\varepsilon_* < 2^{\aleph_0}$ we have: if $\alpha < \beta$ are
from ${\cU}$ then $\varepsilon_{\alpha,\beta} = \varepsilon_*$.  So
if $\alpha < \beta < \gamma$ are from ${\cU}$ then $d_\alpha \in
a_{\gamma,\varepsilon_*} G_{r_\gamma}$ and $d_\beta \in
a_{\gamma,\varepsilon_*} G_{r_\gamma}$, so $d_\alpha =
a_{\gamma,\varepsilon_*} a_1$ and $d_\beta = a_{\gamma,\varepsilon_*} a_2$ for
some $a_1,a_2 \in G_{r_\gamma}$ hence $(d^{-1}_\alpha d_\beta) =
a^{-1}_1 a^{-1}_{\gamma,\varepsilon_*} a_{\gamma,\varepsilon_*} a_2 = a^{-1}_1
a_2 \in G_{r_\gamma}$.]
\mn
\begin{enumerate}
\item[$\circledast_5$]   Without loss of generality
for $\alpha,\beta < \theta$ we have
$d_\alpha \in G_{r_\beta} \Leftrightarrow \alpha \ne \beta$.
\end{enumerate}
\mn
[Why?  Let ${\cU}$ be as in $\circledast_4$.
Without loss of generality otp$({\cU}) = \theta$ and let 
$\langle \alpha_\varepsilon:\varepsilon < \theta\rangle$
list ${\cU}$ in increasing order; let 
$d'_\varepsilon = d^{-1}_{\alpha_{2 \varepsilon}}
d_{\alpha_{2 \varepsilon +1}}$ and let $r'_\varepsilon =
r'_\varepsilon(x,\bar c'_\varepsilon,\bar
b) = r_{\alpha_{2 \varepsilon}}(x,\bar c_{\alpha_{2 \varepsilon}},\bar b)$.

So if $\zeta < \varepsilon < \theta$ then $d^{-1}_{\alpha_{2
\varepsilon}} d_{\alpha_{2 \varepsilon +1}} \in G[r_{\alpha_{2 \zeta}}]$
by $\circledast_1(c)$ hence $d'_\varepsilon = d^{-1}_{\alpha_{2
\varepsilon}} d_{\alpha_{2 \varepsilon +1}} \in G_{r_{\alpha_{2
\zeta}}} = G_{r'_\zeta}$.  Also if $\varepsilon < \zeta < \theta$ then
$d'_\varepsilon = d^{-1}_{\alpha_{2 \varepsilon}} 
d_{\alpha_{2 \varepsilon +1}} \in G_{r_{2 \zeta}} = G_{r'_\zeta}$
by $\circledast_4$.

Also, if $\varepsilon = \zeta$ then $d_{\alpha_{2 \varepsilon +1}} \in
G_{r_{\alpha_{2 \varepsilon}}}$ by $\circledast_1(c)$ and 
$d_{\alpha_{2 \varepsilon}} \notin G_{r_{\alpha_{2 \varepsilon}}}$ by
$\circledast_1(d)$  hence $d'_\varepsilon =
d^{-1}_{\alpha_{2 \varepsilon}} d_{\alpha_{2 \varepsilon +1}} \notin
G_{r_{\alpha_{2 \varepsilon}}} = G_{r'_{\alpha_\varepsilon}}$.  Of
course, $d_{\alpha_{2 \varepsilon}},d_{\alpha_{2 \varepsilon+1}} \in 
G_{\cup\{\bar c_\beta:\beta < \alpha_{2 \varepsilon}\} \cup \bar b}$ hence
$d'_\varepsilon = d^{-1}_{\alpha_{2 \varepsilon}} d_{\alpha_{2
\varepsilon +1}} \in G_{\cup\{\bar c'_\zeta:\zeta < \varepsilon\} \cup
\bar b}$.

Moreover $d'_\varepsilon \notin \psi_1({\gC},\bar c'_\varepsilon,\bar
b)$ as otherwise we recall that $d_{\alpha_{2 \varepsilon}} =
d_{\alpha_{2 \varepsilon +1}}(d'_\varepsilon)^{-1}$ and
$d_{\alpha_{2\varepsilon+1}} \in r_{\alpha_{2\varepsilon}}({\gC},
\bar c_{\alpha_{2\varepsilon}},\bar b) \subseteq \psi_1({\gC},
\bar c_{\alpha_{2\varepsilon}},\bar b)$ and, we are now assuming 
$d'_\varepsilon \in
\psi_1({\gC},\bar c'_\varepsilon,\bar b) = \psi_1({\gC},\bar
c_{\alpha_{2\varepsilon}},\bar b)$ together by $\circledast_1(f)$ we have
$d_{\alpha_{2\varepsilon}} \in \psi_0({\gC},\bar
c_{\alpha_{2\varepsilon}},\bar b)$, contradiction to
$\circledast_1(d)$.  So letting $\psi'_n = \psi_{n+1}$, clearly
renaming we are done as we shall not use $\alpha \ge \theta$.] 

Now by induction on $\varepsilon < \kappa$ we choose
$A_\varepsilon,\bar b_\varepsilon,\langle d_{\alpha,\varepsilon}:\alpha <
\theta\rangle$ from the model $M$ such that:
\mn
\begin{enumerate}
\item[$\circledast_6$]   $(a) \quad \bar b_\varepsilon \char 94 \langle 
d_{\alpha,\varepsilon}:\alpha < \theta\rangle$ 
realizes tp$(\bar b \char 94 \langle d_\alpha:\alpha <
\theta\rangle,A_\varepsilon)$
\sn
\item[${{}}$]   $(b) \quad A_\varepsilon = 
\cup\{\bar c_\alpha:\alpha < \theta\} \cup \bigcup\{\bar b_\zeta,
\bar d_{\alpha,\zeta}:\alpha < \theta$ and $\zeta < \varepsilon\} \cup A_*$.
\end{enumerate}
\mn
[Why possible?  Because $M$ is $\kappa$-saturated, see Definition \ref{dt.28}.]

For $\alpha < \theta,\varepsilon < \kappa$, let $r_{\alpha,\varepsilon} :=
\{\psi_n(x,\bar c_\alpha,\bar b_\varepsilon):n < \omega\}$, so 
$G_{r_{\alpha,\varepsilon}}$ is a
subgroup of $G$ of bounded index (even $\le 2^{\aleph_0}$).
Now for $\alpha,\varepsilon < \theta$ clearly
$d_\alpha \in \cap\{G_{r_{\beta,\zeta}}:\beta < \theta,\zeta <
\varepsilon\}$ by clause $\circledast_1(e)$ (as $r_{\beta,\zeta} 
\in \mathbf R_M$).  Hence by the choice of $d_{\alpha,\varepsilon}$ as
realizing tp$(d_\alpha,A_\varepsilon)$, see $\circledast_6(a)$, as
$A_\varepsilon \supseteq \cup\{\text{Dom}(r_{\beta,\zeta}):\beta <
\theta,\zeta < \varepsilon\}$, clearly
\mn
\begin{enumerate}
\item[$\circledast_7$]  $d_{\alpha,\varepsilon} \in
\cap\{G_{r_{\beta,\zeta}}:\beta < \theta$ and $\zeta < \varepsilon\}$
by $\circledast_6(b)$.
\end{enumerate}
\mn
But $d_\alpha \in \cap\{G_{r_\beta}:\beta < \theta$ and $\beta \ne \alpha\}$
by $\circledast_5$ and $\bar b_\varepsilon \char 94 \langle
d_{\beta,\varepsilon}:\beta < \theta\rangle$ realizes tp$(\bar b \char
94 \langle d_\beta:\beta < \theta\rangle,A_\varepsilon)$ by $\circledast_6$, so
\mn
\begin{enumerate}
\item[$\circledast_8$]   $d_{\alpha,\varepsilon} \in
\cap\{G_{r_{\beta,\zeta}}:\beta < \theta,\beta \ne \alpha$ and $\zeta =
\varepsilon\}$.
\end{enumerate}
\mn
Also by $\circledast_6(a) + \circledast_1(d) + \circledast_2$ we have
\mn
\begin{enumerate}
\item[$\circledast_9$]   $d_{\alpha,\varepsilon} \notin
G_{r_{\alpha,\varepsilon}}$ moreover ${\gC} \models \neg
\psi_0(d_{\alpha,\varepsilon},\bar c_\alpha,\bar b_\varepsilon)$.
\end{enumerate}
\mn
Also by $\circledast_1(f) + \circledast_2$ we have
\mn
\begin{enumerate}
\item[$\circledast_{10}$]   if $d_1,d_2 \in \psi_{n+1}({\gC},
\bar c_\alpha,\bar b'_\varepsilon)$ then
$d_1 d_2,d^{-1}_1,d_1 d^{-1}_2 \in \psi_n({\gC},
\bar c_\alpha,\bar b'_\varepsilon)$.
\end{enumerate}
\mn
Now forget $M$ but retain
\mn
\begin{enumerate}
\item[$(*)_1$]   $\psi_n(n < \omega),b_\varepsilon,\bar
c_\varepsilon,r_{\alpha,\varepsilon},d_{\alpha,\varepsilon},r_\alpha(\alpha <
\theta,\varepsilon < \kappa)$ satisfy $\circledast_1(b)-(f),
\circledast_7,\circledast_8,\circledast_9,\circledast_{10}$.
\end{enumerate}
\mn
Now by Ramsey theorem and compactness, \wilog
\mn
\begin{enumerate}
\item[$(*)_2$]   $\left<\langle \bar c_\alpha:\alpha < \theta\rangle 
\char 94 \bar b_\varepsilon \char 94 
\langle d_{\alpha,\varepsilon}:\alpha < \theta\rangle:
\varepsilon < \kappa\right>$ is an indiscernible sequence over $A_*$.
\end{enumerate}
\mn
Let $d'_{\alpha,\varepsilon} = d^{-1}_{\alpha,2 \varepsilon}
d_{\alpha,2 \varepsilon +1}$ and $\bar b'_\varepsilon = 
\bar b_{2 \varepsilon}$ for $\varepsilon < \kappa,\alpha < \theta$ and let
$r'_{\alpha,\varepsilon}(x) = r(x,\bar c_\alpha,b'_\varepsilon) 
= r_{\alpha,2 \varepsilon}$.  

Now
\mn
\begin{enumerate}
\item[$(*)_3$]   $d'_{\alpha,\varepsilon} \in r({\gC},
\bar c_\beta,\bar b'_\zeta) \cap G$ \underline{iff} $(\alpha,\varepsilon) \ne
(\beta,\zeta)$.
\end{enumerate}
\mn
[Why?  First assume $\varepsilon > \zeta$ then by $(*)_1 +
\circledast_7$ we have
$d_{\alpha,2 \varepsilon},d_{\alpha,2 \varepsilon +1} \in 
G_{r_{\beta,2 \zeta}}$ hence $d'_{\alpha,\varepsilon} = 
d^{-1}_{\alpha,2 \varepsilon} d_{\alpha,2 \varepsilon +1} \in 
G_{r_{\beta,2\zeta}} = G_{r'_{\beta,\zeta}}$.

Second, assume $\varepsilon < \zeta$, then by the indiscernibility,
i.e. $(*)_2$ easily $d_{\alpha,2 \varepsilon} G_{r_{\beta,2 \zeta}} = 
d_{\alpha,2 \varepsilon +1} G_{r_{\beta,2 \zeta}}$ hence 
$d_{\alpha,2 \varepsilon} G_{r'_{\beta,\zeta}} = d_{\alpha,2
\varepsilon +1} G_{r'_{\beta,\zeta}}$ so $d'_{\alpha,\varepsilon} = 
d^{-1}_{\alpha,2 \varepsilon} d_{\alpha,2 \varepsilon +1} \in
G_{r'_{\beta,\zeta}}$ as required.          

Third, assume $\varepsilon = \zeta,\alpha \ne \beta$, then
we have $d'_{\alpha,\varepsilon} \in r({\gC},
\bar c_\beta,b'_\zeta)$ because: $d_{\alpha,2 \varepsilon}
\in r({\gC},\bar c_\beta,\bar b'_\zeta)$ as $\bar b'_\zeta = b_{2
\zeta} = b_{2 \varepsilon}$ in the present case 
and as $\alpha \ne\beta$ using $(*)_1 + \circledast_8$ and
$d_{\alpha,2 \varepsilon +1} \in r({\gC},\bar c_\beta,\bar
b'_\zeta)$ as $\bar b'_\zeta = \bar b_{2 \zeta} = 
\bar b_{2 \varepsilon}$ and as $2 \zeta = 2 \varepsilon < 2 \varepsilon +1$
in the present case , by $\circledast_7 + (*)_1$; of course, $d_{\alpha,2
\varepsilon},d_{\alpha,2 \varepsilon+1} \in G$ hence
$d'_{\alpha,\varepsilon} = d^{-1}_{\alpha,2 \varepsilon} d_{\alpha,2
\varepsilon +1} \in r({\gC},\bar c_\beta,\bar b'_\zeta) \cap G$.

Fourth, assume $\varepsilon = \zeta,\alpha = \beta$.  So 
by $\circledast_9 + (*)_1$, we know that $d_{\alpha,2 \varepsilon} \notin
r({\gC},\bar c_\alpha,\bar b_{2 \varepsilon})$ which means
$d_{\alpha,2 \varepsilon} \notin r({\gC},\bar c_\alpha,\bar
b'_\varepsilon)$.  By $\circledast_7 + (*)_1$ we know that
$d_{\alpha,2 \varepsilon +1} \in r({\gC},\bar c_\alpha,
\bar b_{2 \varepsilon}) = r({\gC},\bar c_\alpha,\bar
b'_\varepsilon)$ and of course $d_{\alpha,2 \varepsilon},d_{\alpha 2
\varepsilon+1} \in G$.
Putting together the last two sentences and the choice of
$d'_{\alpha,\varepsilon}$, as $r({\gC},\bar c_\alpha,
\bar b'_\varepsilon) \cap G$ is a subgroup of $G$ we have
$d'_{\alpha,\varepsilon} = (d_{\alpha,2 \varepsilon})^{-1} d_{\alpha,2
\varepsilon +1} \notin r({\gC},\bar c_\alpha,\bar b'_\varepsilon)
\cap G$ as required in $(*)_3$.]
\mn
\begin{enumerate}
\item[$(*)_4$]   $d'_{\alpha,\varepsilon} \in \psi_1({\gC},\bar
c_\beta,\bar b'_\zeta)$ iff $(\alpha,\varepsilon) \ne (\beta,\zeta)$.
\end{enumerate}
\mn
[Why?  The ``if" direction holds by $(*)_3$ because $r(x,\bar
c_\beta,\bar b'_\zeta) = \{\psi_n(x,\bar c_\beta,\bar b_\zeta):n <
\omega\}$.  For the other direction assume $(\alpha,\varepsilon) =
(\beta,\zeta)$.   As in the proof of $(*)_3$ we have 
$d_{\alpha, 2\varepsilon} \notin \psi_0({\gC},\bar
c_\alpha,\bar b'_\varepsilon)$.  Also $d_{\alpha,2 \varepsilon} =
d_{\alpha,2 \varepsilon+1}(d'_{\alpha,\varepsilon})^{-1}$ so
$d_{\alpha,2 \varepsilon +1}(d'_{\alpha,\varepsilon})^{-1} \notin
\psi_0({\gC},\bar c_\alpha,\bar b_{2 \varepsilon})$.  Also as in
the proof of $(*)_3$ we have $d_{\alpha,2 \varepsilon +1} \in r({\gC},
\bar c_\alpha,\bar b'_\varepsilon)$ hence $d_{\alpha,2\varepsilon
+1} \in \psi_1({\gC},\bar c_\alpha,\bar b'_\alpha)$.  By the last
two sentences and $\circledast_{10} + (*)_1$
we have $d'_{\alpha,\varepsilon} \in \psi_1({\gC},
\bar c_\alpha,\bar b'_\varepsilon) \Rightarrow
d_{\alpha,2\varepsilon} \in \psi_0({\gC},\bar c_\alpha,\bar
b'_\varepsilon)$ but we already note the conclusion fails so
$d'_{\alpha,\varepsilon} \notin \psi_1({\gC},\bar c_\alpha,
\bar b'_\varepsilon)$.

So we are done proving $(*)_4$.]
\mn
\begin{enumerate}
\item[$(*)_5$]   if $u \subseteq \theta \times \kappa$ then some
$d \in {\gC}$ realizes $\{\psi_3(x,\bar c_\alpha,\bar
b'_\zeta)^{\text{if}((\alpha,\zeta) \in u)}:\alpha <
\theta,\varepsilon < \kappa\}$.
\end{enumerate}
\mn
[Why?  By saturation \wilog \, $u$ is co-finite, let $\langle
(\alpha(\ell),\varepsilon(\ell)):\ell < k\rangle$ list $\theta \times
\kappa \backslash u$ with no repetitions and let $d := 
d'_{\alpha(0),\varepsilon(0)} d'_{\alpha(1),\varepsilon(1)} \ldots
d'_{\alpha(k-1),\varepsilon(k-1)}$.

On the one hand by $(*)_3$ clearly 

\begin{equation*}
\begin{array}{clcr}
(\alpha,\varepsilon) \in u &\Rightarrow
(\alpha,\varepsilon) \in \theta \times \kappa \backslash
\{(\alpha(\ell),\varepsilon(\ell)):\ell < k\}\\
  &\Rightarrow \{d'_{\alpha(\ell),\varepsilon(\ell)}:\ell < k\} 
\subseteq r({\gC},\bar c_\alpha,\bar b'_\varepsilon) \cap G \\
  &\Rightarrow d \in r({\gC},\bar c_\alpha,\bar b_\varepsilon)
\Rightarrow {\gC} \models \psi_3[d,\bar c_\alpha,\bar b'_\varepsilon].
\end{array}
\end{equation*}

\mn
On the other hand if $\ell < k$ then let $e_1 =
d'_{\alpha(0),\varepsilon(0)} \ldots
d'_{\alpha(\ell-1),\varepsilon(\ell-1)}$ and let $e_2 =
d'_{\alpha(\ell+1),\varepsilon(\ell+1)} \ldots
d'_{\alpha(k-1),\varepsilon(k-1)}$ so $d =
e_1 d'_{\alpha(\ell),\varepsilon(\ell)} e_2$ hence
$d'_{\alpha(\ell),\varepsilon(\ell)} = e^{-1}_1 de^{-1}_2 $.  As above
$e_1,e_2 \in r({\mathfrak C},\bar c_{\alpha(\ell)},
\bar b'_{\varepsilon(\ell)}) \cap G$ hence 
$e^{-1}_1,e^{-1}_2 \in \psi_3({\gC},\bar c_{\alpha(\ell)},
\bar b_{\varepsilon(\ell)}) \cap G$.  
As $d \in G$ and $e^{-1} d \in G$, by $\circledast_{10}$ we get

\begin{equation*}
\begin{array}{clcr}
d \in \psi_3({\gC},\bar c_{\alpha(\ell)},
\bar b'_{\varepsilon(\ell)}) &\Rightarrow
e^{-1}_1 d \in \psi_2({\gC},\bar c_{\alpha(\ell)},
\bar b'_{\varepsilon(\ell)}) \\
 &\Rightarrow e^{-1}_1 d e^{-1}_2 \in \psi_1({\gC},\bar
c_{\alpha(\ell)},\bar b'_{\varepsilon(\ell)}) \\
  &\Rightarrow d'_{\alpha(\ell),\varepsilon(\ell)} \in 
\psi_1({\gC},\bar c_{\alpha(\ell)},\bar b'_{\varepsilon(\ell)}).
\end{array}
\end{equation*}

\mn
But this contradicts $(*)_4$.

Now $(*)_5$ gives ``$\psi_3(x,\bar z,\bar y)$ witness $T$ is
2-independent" so we are done.    
\end{PROOF}

\begin{claim}
\label{dt.32}
If $G$ is Abelian, then \ref{dt.35}
can be proved also replacing $q_B({\gC})$ by $\cap\{G':G'$ is a subgroup of
$G$ of bounded index preserved by automorphisms of ${\gC}$ over $A
\cup A_*\}$.
\end{claim}

\begin{proof}  We shall prove this elsewhere.
\end{proof}

\begin{discussion}
\label{dt.49}
1) Is $1/2$-dependence preserved by weak
expansions (as in \cite[\S1]{Sh:783})?  Of course not, as if $M$ is a
model of $T_2$ from \ref{dt.23} then any $Y \subseteq 
\prod\limits_{\ell <k} P^M_i$ is definable in such an expansion, and easily
for some such $Y$ we can interpret number theory (as number theory is
interpretable in some bi-partite graph).

\noindent
2) Is the following interesting?   I think yes!  It seems
that we can prove the $k$-dimensional
version of \ref{dt.35} for $1/k$-dependent $T$, i.e. for
$k=1$ it should give \cite{Sh:876}, for $k=2$ it should give
\ref{dt.35}.  E.g. think of having $|T| < \lambda_0 < \ldots <
\lambda_k,\lambda_{\ell +2} = (\lambda_{\ell +1})^{\lambda_\ell}$ and
we choose $M_\ell \prec {\gC}$ of cardinality $\lambda_\ell$
closed enough by downward induction on $\ell$, i.e. we get 
a ${\cP}^-(k)$-diagram.  We shall try to deal with this elsewhere.
\end{discussion}
\newpage

\bibliographystyle{amsalpha}
\bibliography{shlhetal}

\end{document}